\documentclass[11pt,
]{article}

\usepackage{amssymb,amsfonts,amstext,amsmath,amsthm,latexsym,mathrsfs}
\usepackage{mathbbol}

\usepackage{makeidx,mparhack}

\makeindex

\makeatletter
\if@twoside%
\else\fi%
\makeatother  

\newcommand{\nek}{\newcommand}
\nek{\renek}{\renewcommand}
\nek{\vyk} [1] {}

\nek{\ubf}{\fontseries{b}\selectfont}
\nek{\bfit}{\bfseries\itshape}
\nek{\bftt}{\ttfamily\bfseries\upshape\selectfont}

\nek{\parf}{\subsection}

\nek{\bsec} {\parf}
\nek{\punk}{\subsubsection}

\renek{\thesubsection}{\thesection\hspace*{0.1ex}\alph{subsection}}

\theoremstyle{plain}

\newtheorem{theorem}             {Theorem} [subsection]
\newtheorem{corollary}  [theorem]{Corollary}
\newtheorem{prop}  [theorem]{Proposition}
\newtheorem{lemma}      [theorem]{Lemma}

\theoremstyle{definition}

\newtheorem{definition} [theorem]{Definition}
\newtheorem{prim}       [theorem]{Example}
\newtheorem{zam}        [theorem]{Remark}
\newtheorem*{prF}   {Proof}    
\newtheorem*{ack}   {Acknowledgement}    
\newtheorem{ggi}   [theorem] {Blanket Assumption}  

\nek{\thsp}{\hspace{0.1ex}}
\nek{\back}{\begin{ack}}
\nek{\eack}{\end{ack}}
\nek{\bpro}{\begin{prop}}
\nek{\epro}{\end{prop}}
\nek{\bcor}{\begin{corollary}}
\nek{\ecor}{\end{corollary}}
\nek{\bdf} {\begin{definition}}
\nek{\edf} {\qed\end{definition}}
\nek{\eDf} {\end{definition}}
\nek{\bgg} {\begin{ggi}}
\nek{\egg} {\qed\end{ggi}}
\nek{\ble} {\begin{lemma}}
\nek{\ele} {\end{lemma}}
\nek{\bpri}{\begin{prim}}
\nek{\epri}{\qed\end{prim}}
\nek{\bte} {\begin{theorem}}
\nek{\ete} {\end{theorem}}
\nek{\bre}{\begin{zam}}
\nek{\ere}{\qed\end{zam}}
\nek{\bpf} {\begin{prF}} 
\nek{\epf} {\qed\end{prF}} 
\nek{\epF}[1]{\hfill\hbox{$\square$\ ({\small#1})}\end{prF}}

\nek{\qeD}[1]{\hfill\hbox{$\square$\ ({\small#1})}}

\nek{\ben}{\begin{enumerate}}
\nek{\een}{\end{enumerate}}
\nek{\bit}{\begin{itemize}}
\nek{\eit}{\end{itemize}}
\nek{\bay}{\begin{array}}
\nek{\eay}{\end{array}}
\nek{\bce}{\begin{center}}
\nek{\ece}{\end{center}}
\nek{\lh}  {\mathop{\tt lh}}
\nek{\dom} {\mathop{\tt dom}}
\nek{\ran} {\mathop{\tt ran}}
\nek{\HC}  {{\text{\rm{HC}}}}
\nek{\hc}{\HC}
\nek{\rL} {{\mathbf L}}
\nek{\rV} {{\mathbf V}}
\nek{\ZFC} {\text{\ubf ZFC}}
\nek{\zfc} {\ZFC}
\nek{\al}  {\alpha}
\nek{\ga}  {\gamma}
\nek{\Ga}  {\Gamma}
\nek{\da}  {\delta}
\nek{\Da}  {\Delta}
\nek{\kpa} {\kappa}
\nek{\la}  {\lambda}
\nek{\ve}  {\varepsilon}
\nek{\vpi} {\varphi}
\nek{\sg}  {\sigma}
\nek{\Sg}  {\Sigma}
\nek{\om}  {\omega}
\nek{\Om}  {\Omega}
\nek{\lom} {^{<\om}} 
\nek{\omi} {\om_1}
\nek{\za}  {\zeta}
\nek{\fs}[2]{{\mathbf\Sigma}^{#1}_{#2}}
\nek{\fp}[2]{{\mathbf\Pi}^{#1}_{#2}}
\nek{\fd}[2]{{\mathbf\Delta}^{#1}_{#2}}
\nek{\is}[2]{\varSigma^{#1}_{#2}}
\nek{\ip}[2]{\varPi^{#1}_{#2}}
\nek{\id}[2]{\varDelta^{#1}_{#2}}
\nek{\BBB}{\hspace{0.1ex}}
\nek{\dP}{{\BBB{\mathbb P}\BBB}}
\nek{\dS}{{\BBB{\mathbb S}\BBB}}
\nek{\pu}  {\varnothing}
\nek{\sq}  {\subseteq}
\nek{\qs}  {\supseteq}
\nek{\su}  {\subset}
\nek{\eqv} {\mathbin{\,\Longleftrightarrow\,}}
\nek{\imp} {\mathbin{\,\Longrightarrow\,}}
\nek{\mpi} {\mathbin{\,\Longleftarrow\,}}
\nek{\lra} {\longrightarrow} 
\nek{\we}  {{\mathbin{\hspace*{0.2ex}^\wedge}}}
\nek{\sus} {{\exists\,}}
\nek{\kaz} {{\forall\,}}
\nek{\ti}  {\times}
\nek{\dm}  {$$}
\nek{\iy}  {\infty}
\nek{\abs}[1]{|#1|}

\nek{\nin}{\notin}
\nek{\res} {{\hspace*{0.1ex}\restriction\hspace*{0.3ex}}}
\nek{\ang} [1] {\langle #1\rangle}
\nek{\ans} [1] {\{\hspace{0.1ex}#1\hspace{0.1ex}\}}
\nek{\dd}[1]{$\mtho\hspace{0.2ex}{#1}$-\hspace{0.0ex}}
\nek{\itla}{\item\label}
\nek{\mtho}{\mathsurround=0mm}
\nek{\msur}{\hspace*{-1\mathsurround}}
\nek{\dsur}{\hspace{-0.3\mathsurround}}
\nek{\hsur}{\hspace{-0.5\mathsurround}}
\nek{\noi}{\noindent}
\nek{\vom}{\vspace{1mm}}
\nek{\vtm}{\vspace{2mm}}
\nek{\bez}{\smallsetminus}
\nek{\rit} [1] {{\it#1\/}}

\nek{\ens} [2] {\ans{{#1\hspace{0.5ex}{:}}\zz\hspace{0.5ex}#2}}
\nek{\zz} {\linebreak[0]} 
\nek{\goth}{\mathfrak}
\nek{\skri}{\mathscr}

\nek{\cP}{\skri P}

\nek{\yo} {,\linebreak[0]}
\nek{\yi} {\hspace{\mathsurround},\linebreak[0]\hspace{\mathsurround}}
\nek{\yd} {\hspace{\mathsurround},\linebreak[0]\:}
\nek{\yt} {\hspace*{\mathsurround}\text{,}\linebreak[0]\;}

\nek{\eg} {\text{\sl e.\resp g.}}

\nek{\iesp}{\hspace{0.3ex}}
\nek{\resp}{\hspace{0.25ex}}

\nek{\itsep}{\itemsep=0.25ex plus 0.1ex minus 0.1ex}

\nek{\tenu}[1]{
\def\theenumi{#1}
\itsep}

\nek{\nenu}{\tenu{{\rm(\arabic{enumi})}}}
\nek{\renu}{\itsep\tenu{{\rm(\roman{enumi})}}\itsep}

\nek{\atc} {\addtocounter{enumi}1}

\nek{\stk} [2] {\ang{#1\hspace{0.3ex};\hspace{0.1ex}#2}}
\nek{\sis} [2] {\ans{#1}_{#2}}

\nek{\cM} {\goth M}

\nek{\lam} [1]
{\label{#1}\hspace*{-3pt}\imar{#1}%
}%
\nek{\las} [1]
{\label{#1}\imae{#1}}%

\nek{\imar}[1]{\marginpar[
\flushright\footnotesize%
$\mtho\longrightarrow$\\%
\vspace{-1ex}{#1}\vspace*{1ex}]%
{
\flushleft\footnotesize%
$\mtho\longleftarrow$\\%
\vspace{-1ex}{#1}\vspace*{1ex}}%
}%

\nek{\imae}[1]{\marginpar[
\flushright\footnotesize\vspace{-4ex}%
$\mtho\longrightarrow$\\%
\vspace{-1ex}{#1}$\mtho$\vspace*{2ex}]%
{
\flushleft\footnotesize\vspace{-4ex}%
$\mtho\longleftarrow$\\%
\vspace{-1ex}{#1}$\mtho$\vspace*{2ex}}
}%

\nek{\pws} [1] {\cP(#1)}

\nek{\snos} [1] {\,\footnote{\ #1}}
\nek{\snom}   {\,\footnotemark}
\nek{\snot} [1] {\footnotetext{\ #1}}

\nek{\onto} {\overset{{\text{onto}}}{\longrightarrow}}

\nek{\gp} {\mathbf p}

\nek{\limp} {\mathrel{\,\imp\,}}

\nek{\spe} [1] {\text{\bf SS}(#1)}
\nek{\spf} [1] {\text{\bf SS}(#1)}
\nek{\spg} [1] {\text{\bf SS}^{<\om}(#1)}

\nek{\vys} [1] {\text{\tt hgt}(#1)} 

\nek{\cle} {\preccurlyeq}
\nek{\cls} {\prec}

\renek{\refname} {{\large\bf References}}

\nek{\dphi} {\mathbb\Phi}

\nek{\dU}{\mathbb U}
\nek{\uu}{^{\text{\tt fu}}}

\nek{\dn} {2^\om}

\nek{\dox}{{\namy{\boldsymbol x}}}

\nek{\rc} {\mathbf c}
\nek{\rd} {\mathbf d}
\nek{\rpi} {\dox}

\nek{\plo} {\dP\lom}
\nek{\uplo} {(\dP\cup\dU)\lom}

\renek{\thesubsection}{\arabic{subsection}}

\nek{\cll} {\mathrel{{\cls}\hspace*{-0.9ex}{\cls}}}

\nek{\vpj} [2] {\vpi^{#1}_{#2}}
\nek{\Phj} [2] {\Phi^{#1}_{#2}}

\nek{\od} {\text{OD}}

\nek{\Eo} {\mathrel{{\text{\sf E}}_0}}

\nek{\zf} {\text{\ubf ZF}}

\nek{\uf} [2] {{\boldsymbol U}^\dphi_{#1}(#2)}
\nek{\ufi} [1] {{\boldsymbol U}^\dphi_{#1}}

\nek{\tx} [2] {{\boldsymbol T}^\dphi_{#1}(#2)}
\nek{\txi} [1] {{\boldsymbol T}^\dphi_{#1}}

\nek{\lel} {\leq_{\rL}}

\nek{\kc} [2] {C_{#1}^{#2}}
\nek{\xc} [2] {C'_{#1#2}}
\nek{\kd} [2] {D_{#1#2}}
\nek{\kcp} [3] {C_{#1#2}^{#3}}
\nek{\kr} [3] {\jrho_{#1#2}^{#3}}
\nek{\kR} [2] {R^{#1}_{#2}}

\nek{\roo} [1] {\text{\tt stem}(#1)}

\nek{\pes} {\text{\bf ST}}
\nek{\app}{\boldsymbol\cdot}

\nek{\ssptf} {{\ubf STF}}
\nek{\sif} {\ssptf}  

\nek{\dplo} {\dd\dP}

\nek{\ek} [2] {[#1]_{#2}}
\nek{\eke} [1] {[#1]_{\rE}}
\nek{\eko} [1] {[#1]_{\Eo}}
\nek{\rE}{\mathrel{\mathsf E}}

\nek{\mus} {multisystem}

\nek{\bse} {2^{<\om}}

\nek{\La} {\Lambda}
\nek{\leqv} {\,\eqv\,}
\nek{\qand} {\quad\text{and}\quad}
\nek{\sqf} {\sq^{\text{\tt fin}}}

\nek{\namy} [1] 
{\overset{\text{\mtho$\hspace*{0.5ex}_\text{\Large\bf.}$}}{#1}}

\nek{\cD} {\mathscr D}

\nek{\ka}{\kappa}
\nek{\vt}{\vartheta}
\nek{\ret} [1] {\res_{\hspace*{0.05ex}#1}}
\nek{\pet} {\text{\ubf PT}}

\renek{\imar}[1] {}
\renek{\imae}[1] {}

\begin{document}

\title
{A definable $\Eo$ class containing no definable elements
}

\author 
{
Vladimir~Kanovei\thanks{IITP RAS and MIIT,
  Moscow, Russia, \ {\tt kanovei@googlemail.com} --- contact author. 
Partial support of 
RFFI grant 13-01-00006 acknowledged.}  
\and
Vassily~Lyubetsky\thanks{IITP RAS,
  Moscow, Russia, \ {\tt lyubetsk@iitp.ru} 
}}

\date 
{\today}

\maketitle

\begin{abstract}
A generic extension $\rL[x]$ of $\rL$ by a real $x$ is defined, 
in which the $\mathsf E_0$-class of $x$ is a $\ip12$ set 
containing no ordinal-definable reals.
\end{abstract}

\parf{Introduction}
\las{cha1}

It is known that the existence of a non-empty $\od$ 
(ordinal-definable) set of reals $X$ with no \od\  
element is consistent with $\ZFC$; the set of all 
non-constructible reals gives an example in many generic 
models including \eg\ the Solovay model or the extension 
of $\rL$, the constructible universe, by a Cohen real.
\begin{quotation}
\noi
\rit{Can such a set\/ $X$ be countable}? 
That is, is it consistent with $\ZFC$ that there is a countable 
$\od$ (or outright definable by a precise set-theoretic formula) 
set of reals $X$ containing no \od\ element?
\end{quotation}
This question was initiated and 
discussed at the \rit{Mathoverflow} website\snos
{\label{snos1}
A question about ordinal definable real numbers. 
\rit{Mathoverflow}, March 09, 2010. 
{\tt http://mathoverflow.net/questions/17608}. 
}
and at FOM\snos
{\label{snos2}%
Ali Enayat. Ordinal definable numbers. FOM Jul 23, 2010.
{\tt http://cs.nyu.edu/pipermail/fom/2010-July/014944.html}}
.
In particular Ali Enayat (Footnote~\ref{snos2}) conjectured that 
the problem can be solved by the   
finite-support countable product $\plo$ 
of the Jensen ``minimal $\ip12$ real 
singleton forcing'' $\dP$ defined in \cite{jenmin} 
(see also Section 28A of \cite{jechmill}).
Enayat proved that a symmetric part of the \dd\plo generic 
extension of $\rL$ definitely yields a model of $\zf$ 
(not a model of $\ZFC$!) 
in which there is a Dedekind-finite infinite \od\ set of 
reals with no \od\ elements --- namely the set of all 
reals \dd\dP generic over $\rL$. 
In fact \dd\plo generic extensions of $\rL$ and their symmetric 
submodels were considered in \cite{ena} (Theorem 3.3) with 
respect to some other questions. 

Following the mentioned conjecture, we proved in \cite{kl:cds} 
that indeed, in a \dd\plo generic extension of $\rL$, 
the set of all reals \dd\dP generic over $\rL$ 
is a countable $\ip12$ set with no OD elements. 
The $\ip12$ definability is definitely the best one can get 
in this context since it easily follows from the $\ip11$ 
uniformisation theorem that any non-empty $\is12$ set 
of reals definitely contains a $\id12$ element.

Jindra Zapletal\snos
{\label{snos3}%
Personal communication, Jul 31/Aug 01, 2014.} 
informed us that there is a totally different model of $\ZFC$ 
with an $\od$ \dd\Eo class\snos
{Recall that if $x,y\in\om^\om$ then $x\Eo y$ iff $x(n)=y(n)$ 
for all but finite $n$.}
$X$ containing no \od\ elements.
The construction of such a model, not yet published, 
but described to us in a brief communication, 
involves a combination of several forcing notions and 
some modern ideas in descriptive set theory, like models of 
the form $\rV[x]_{\rE}$ for ${\rE}={\Eo}$,  
recently presented in \cite{ksz}; it also does not look to 
yield $X$ being analytically definable, let alone $\ip12$.

We prove the next theorem in this paper$:$

\bte
\lam{mt}
It is true in a suitable generic extension\/ $\rL[x]$ 
of\/ $\rL$, the constructible universe, by a real\/ $x\in\dn$ 
that the\/ \dd\Eo equivalence class\/ $\eko x$ 
(hence a countable set) 
is\/ $\ip12$, but it has no\/ $\od$ elements.
\ete

The forcing $\dP$ we use to prove the theorem is a clone of 
the abovementioned Jensen forcing, 
but defined on the base of the
Silver forcing instead of the Sacks forcing.  
The crucial advantage of Silver's forcing here is that it
leads to  
a Jensen-type forcing naturally closed under the 0-1 flip 
at any digit, so that the corresponding extension contains 
a $\ip12$ \dd\Eo class of generic reals instead of a 
$\ip12$ generic singleton as in \cite{jenmin}.
In fact a bigger family of \rit{\dd\Eo large\/} trees
(perfect trees $T\sq\bse$ such that ${\Eo}\res{[T]}$ 
is not smooth, see \cite[Section 10.9]{kanB}) 
would also work similarly to Silver trees, an by similar
reasons.

\bre
\lam{mtr}
Theorem \ref{mt} also solves another question asked
at the \rit{Mathoverflow} website\snos
{\label{snos4}
A question about definable non-empty sets containing no
definable elements. 
\rit{Mathoverflow}, February 11, 2013, 
{\tt http://mathoverflow.net/questions/121484}. 
}~:
namely, 
\begin{quotation}
\noi
is there an example of a set $S$  
definable in $\ZFC$ and provable in $\ZFC$ to be
countably infinite, while at the same time,
no set definable in $\ZFC$ can be proved in $\ZFC$
to be an element of $S$?
\end{quotation}
To define such an example, let $S$ be defined as
(1) $\eko x$ provided
the set universe is equal to the class $\rL[x]$ as in Theorem~\ref{mt},
and
(2) simply $S=\om$ otherwise.
Suppose towards the contrary that $\ZFC$ proves that the 
real $x$, uniquely defined by a certain fixed formula, 
outright belongs to $S$. 
Then in particular this must be true in case (1), contrary to the 
definition of $S$ via Theorem~\ref{mt}.
\ere

It remains to note that a \rit{finite} \od\ set of reals 
contains only \od\ reals by obvious reasons. 
On the other hand, by a result in \cite{gl} 
there can be two \rit{sets} of reals $X,Y$  
such that the pair $\ans{X,Y}$ is \od\ but 
neither $X$ nor $Y$ is \od.  

\parf{Trees and Silver-type forcing} 
\las{tre}
\label{ptf}

Let $\bse$ be the set of all strings (finite sequences) 
of numbers $0,1$.
If $t\in\bse$ and $i=0,1$ then 
$t\we k$ is the extension of $t$ by $k$. 
If $s,t\in\bse$ then $s\sq t$ means that $t$ extends $s$, while 
$s\su t$ means proper extension. 
If $s\in\bse$ then $\lh s$ is the length of $s$,  
and $2^n=\ens{s\in\bse}{\lh s=n}$ (strings of length $n$).%

Let any $s\in\bse$ {\ubf act} on $\dn$ so that
$(s\app x)(k)=x(k)+s(k)\pmod 2$ whenever $k<\lh s$ and simply
$(s\app x)(k)=x(k)$ otherwise.
If $X\sq\dn$ and $s\in\bse$ then, as usual, let
$s\app X=\ens{s\app x}{x\in X}$.

Similarly if $s\in 2^m,$ $t\in2^n,$ $m\le n$, then define
$s\app t\in 2^n$ so that
$(s\app t)(k)=t(k)+s(k)\pmod 2$ whenever $k<\min\ans{m,n}$ and  
$(s\app t)(k)=t(k)$ whenever $m\le k<n$.
Note that $\lh(s\app t)=\lh t$.
Let $s\app T=\ens{s\app t}{t\in T}$ for $T\sq\bse$.

If $T\sq\bse$ is a tree and $s\in T$ then put 
$T\ret s=\ens{t\in T}{s\sq t\lor t\sq s}$. 

Let $\pet$ be the set of all \rit{perfect} trees 
$\pu\ne T\sq \bse$  
\imar{pet}%
(those with no endpoints and no isolated branches). 
If $T\in\pet$ then there is a largest string $s\in T$ 
such that $T=T\ret s$; it is denoted by $s=\roo T$   
(the {\it stem\/} of $T$); 
we have $s\we 1\in T$ and $s\we 0\in T$ in this case.
If $T\in\pet$ then  
$$ 
[T]=\ens{a\in\dn}{\kaz n\,(a\res n\in T)}\sq\dn
$$ 
is the perfect set of all \rit{paths through $T$}.

Let $\pes$ be the set of all \rit{Silver trees}, that is, those
$T\in\pet$ that is a partition
$\om=u_0\cup u_1\cup u_{0,1}$ such that $u_{0,1}$ is infinite and
if $s\in T$ then\vtm

$-$ \ if $\lh s\in u_0$ then $s\we0\in T$ but $s\we1\nin T$;\vtm

$-$ \ if $\lh s\in u_1$ then $s\we1\in T$ but $s\we0\nin T$;\vtm

$-$ \ if $\lh s\in u_{0,1}$ then $s\we0\in T$ and $s\we1\in T$.\vtm

\noi
By a {\ubf Silver-type forcing}
(\ssptf)
we understand any set 
$\dP\sq\pes$ such that 
\ben
\nenu
\itla{ptf1} 
$\dP$ contains the full tree $\bse$;
\imar{ptf1}

\itla{ptf2}
\label{utp} 
if $u\in T\in\dP$ then $T\ret u\in \dP$.
\imar{ptf2}

\itla{ptf3}
if $T\in\dP$ and $s\in\bse$ then $s\app T\in \dP$.
\imar{ptf3}
\een
Such a set $\dP$ can be considered as a forcing notion 
(if $T\sq T'$ then $T$ is a stronger condition), and then 
it adds a real in $\dn$.

\parf{Splitting construction over a Silver-type forcing} 
\las{spe}

Assume that $\dP\sq\pes$ is a \sif.
The set  $\spe\dP$ of \rit{Silver splitting}
constructions over $\dP$ 
consists of all finite systems of trees of the form 
$\vpi=\sis{T_s}{s\in2^{< n}}$, where 
$n=\vys\vpi<\om$ (the height of $\vpi$), satisfying the
following conditions:
\ben
\nenu
\atc\atc\atc
\itla{spe1}
each tree $T_s=\vpi(s)$ belongs to $\dP$, \ --- \ we let
$r_s=\roo{T_{s}}$;
\imar{spe1}

\itla{spe2}
if $s\we i\in 2^{<n}$ ($i=0,1$) then
$T_{s\we i}\sq T_s\ret{r_{s}\we i}$  --- 
\imar{spe2}
it easily follows that 
$[T_{s\we0}]\cap [T_{s\we1}]=\pu$;

\itla{spe3}
there is an increasing sequence of numbers 
$h(0)<h(1)<\dots<h(n-1)$ such that   
$\lh{r_{s}}=h(k)$ whenever $s\in 2^k$ and $k<n$;  
\imar{spe3} 

\itla{spe4}
if $k<m<n$, $u,v\in 2^{m}$, and $h(k)<j<h(k+1)$ then
$r_{u}(j)= r_{v}(j)$.  
\imar{spe4}

\itla{spe5}
if $m< n$, $u,v\in 2^{m}$, and $t\in\bse$ then
$r_{u}\we t\in T_u\leqv r_{v}\we t\in T_v$.  
\imar{spe5}
\een
The tree $T=\bigcup_{s\in 2^{n-1}}T_s$ belongs to $\pes$ in this case.

Let $\vpi,\psi$ be systems in $\spe\dP$.
Say that 
\bit
\item[$-$]
$\vpi$ \rit{extends} $\psi$, symbolically $\psi\cle\vpi$, if 
$n=\vys\psi\le\vys\vpi$ and $\psi(s)=\vpi(s)$ for 
all $s\in2^{<n}$;

\item[$-$]
\rit{properly extends} $\psi$, 
symbolically $\psi\cls\vpi$, if in 
addition $\vys\psi<\vys\vpi$;

\item[$-$]
\rit{reduces} $\psi$, if $n=\vys\psi=\vys\vpi$, 
$\vpi(s)\sq\psi(s)$ for all $s\in 2^{n-1},$ and 
$\vpi(s)=\psi(s)$ for all $s\in 2^{<n-1}$.
\eit
In other words, the reduction allows to shrink trees in the 
top layer of the system, but keeps intact those in the lower 
layers.

Note that $\vpi=\La$ (the empty system) is the only one 
with $\vys\vpi=0$. 
To get a system $\vpi$ with $\vys\vpi=1$ 
(and then $\dom\vpi=\ans{\La}$) 
put $\vpi(\La)=T$,
where $T\in\pes$.
The following lemma leads to systems of bigger height.

\ble
\lam{lsys}
Assume that\/ $\dP\sq\pes$ is a\/ \sif\ and\/ 
$\vpi=\sis{T_s}{s\in2^{<n}}\in\spe\dP$. 
\ben
\renu
\itla{sys1}
If\/ $s_0\in2^{n-1},$ and\/ $T\in\pes\yt T\sq T_{s_0}$, 
then there is
a system\/ $\vpi'=\sis{T'_s}{s\in2^{<n}}\in\spe\dP$ which
reduces\/ $\vpi$ and satisfies\/ $T_{s_0}=T$.

\itla{sys2}
There is
a system\/ $\vpi'=\sis{T'_s}{s\in2^{<n+1}}\in\spe\dP$ which
properly extends\/ $\vpi$.

\itla{sys3}
If a system\/ $\psi$ properly extends\/ $\vpi$ and 
a system\/ $\psi'$ 
reduces\/ $\psi$ then\/ $\psi'$ properly extends\/ $\vpi$.
\een
\ele
\bpf
By definition all strings $r_s=\roo{T_s}$ with $s\in2^{n-1}$
satisfy $\lh{r_s}=h$ for one and the same $h=h(n-1)$.

\ref{sys1} 
Put $T'_s=\ens{r_s\we t}{r_{s_0}\we t\in T}$ for all 
$s\in 2^{n-1}$, and still $T'_s=T_s$ for $s\in 2^{<n-1}$.
The sets $T'_s$ defined this way 
belong to $\dP$ by \ref{ptf3} of Section \ref{tre}.

\ref{sys2} 
Put 
$T'_{s\we i}=T_s\ret{r_s\we i}$ 
for all 
$s\in2^{n-1}$ and $i=0,1$, and still $T'_s=T_s$ for $s\in 2^{< n}$.
The sets $T'_{s\we i}$ belong to $\dP$ by \ref{ptf2} 
of Section \ref{tre}.
\epf

By the lemma, if $\dP\sq\pes$ is a\/ \sif\ then there is a 
strictly 
\dd\cls increasing sequence $\sis{\vpi_n}{n<\om}$ in 
$\spf\dP$. 
The limit system $\vpi=\bigcup_n\vpi_n=\sis{T_s}{s\in\bse}$
then satisfies conditions \ref{spe1} --- \ref{spe5} 
on the whole domain $\bse$.

\bpro
\lam{issil} 
In this case, the tree\/ $T=\bigcap_n\bigcup_{s\in2^n}T_s$ is 
still a Silver tree in\/ $\pes$ 
{\rm(not necessarily in $\dP$)}, 
and\/ $[T]=\bigcap_n\bigcup_{s\in2^n}[T_s]$.\qed
\epro

Say that a tree $T$ \rit{occurs in\/ $\vpi\in\spf\dP$} if 
$T=\vpi(s)$ for some $s\in 2^{\le\vys\vpi}$.

We define $\spg\dP$, 
{\ubf the finite-support product} of countably many
copies of $\spf\dP$,  
to consist of all infinite sequences  
$\Phi=\sis{\vpi_k}{k\in\om}$, where each $\vpi_k=\Phi(k)$ 
belongs to $\spf\dP$
and the set $\abs \Phi=\ens{k}{\vpi_k\ne\La}$
(the support of $\Phi$) is finite. 
Sequences $\Phi\in\spf\dP$ will be called \rit{\mus s}.  

Say that a tree $T$ \rit{occurs in\/ 
$\Phi=\sis{\vpi_k}{}$} if 
it occurs in some $\vpi_k\yd k\in\abs\Phi$.

Let $\Phi,\Psi$ be \mus s in $\spg\dP$. 
We define that 
\bit
\item[$-$]
$\Phi$ \rit{extends} $\Psi$, symbolically $\Psi\cle\Phi$, if 
$\Psi(k)\cle\Phi(k)$ 
(in $\spf\dP$) for all $k$;

\item[$-$]
$\Psi\cll\Phi$, iff $\abs\Psi\sq\abs\Phi$ 
and $\Psi(k)\cls\Phi(k)$ for all $k\in\abs\Psi$;

\item[$-$]
$\Phi$ \rit{reduces} $\Psi$ iff 
$\Phi(k)$ reduces $\Psi(k)$ for all $k\in\abs\Psi$.
\eit

\bcor
[of Lemma~\ref{lsys}]
\lam{csys}
If\/ $\dP\sq\pes$ is a\/ \sif\ and\/ $\Psi\in\spg\dP$  
then there is a \mus\/ $\Phi\in\spg\dP$ such that\/ 
$\Psi\cll\Phi$.\qed 
\ecor

\parf{Jensen's extension of a Silver-type forcing} 
\las{jex}

Let $\zfc'$ be the subtheory of $\zfc$ including all 
axioms except for the 
power set axiom, plus the axiom saying that $\pws\om$ exists. 
(Then $\omi$ and continual sets like $\pet$ exist as well.)
Let $\cM$ be a countable transitive model of $\ZFC'$. 

Suppose that $\dP\in\cM\yd \dP\sq\pes$ is a \sif.  
Then the sets 
$\spf\dP$ and $\spg\dP$ belong to 
$\cM$, too. 

\bdf
\lam{dPhi}
Consider any \dd\cle increasing sequence 
$\dphi=\sis{\Phi^j}{j<\om}$ of \mus s 
$\Phi^j=\sis{\vpj jk}{k\in\om}\in\spg\dP$, 
\rit{generic over\/ $\cM$} in the sense that it intersects 
every set $D\in\cM\yd D\sq\spg\dP$, dense in $\spg\dP$\snos
{Meaning that for any $\Psi\in\spg\dP$ there is $\Phi\in D$ 
with $\Psi\cle\Phi$.} 
.

Then in particular it intersects every set of the form   
$$
D_k=\ens{\Phi\in\spg\dP}{\kaz k'\le k\: 
(k\le\vys{\Phi(k')}}\,. 
$$
Hence if $k<\om$ then the sequence 
$\sis{\vpj jk}{j<\om}$ of systems $\vpj jk\in\spf\dP$ is 
\rit{eventually strictly increasing}, so that 
$\vpj jk\cls \vpj {j+1}k$ 
for infinitely many indices $j$  
(and $\vpj jk=\vpj {j+1}k$ for other $j$).
Therefore there is a system of trees 
$\sis{\tx ks}{k<\om\land s\in\bse}$ in $\dP$ 
such that $\vpj jk=\sis{\tx ks}{s\in 2^{< h(j,k)}}$, 
where $h(j,k)=\vys{\vpj jk}$.
Then 
$$
\textstyle
\ufi k=\bigcap_n\bigcup_{s\in2^n}\tx ks\qand 
\uf ks=\bigcap_{n\ge \lh s}
\bigcup_{t\in2^n,\:s\sq t}\tx kt
$$   
are trees in $\pes$ (not necessarily in $\dP$) 
by Proposition~\ref{issil} 
for each $k$ and $s\in\bse;$ 
thus $\ufi k=\uf k\La$. 
In fact $\uf ks=\ufi k\cap \tx ks$ by \ref{spe2}.
\edf

\ble
\lam{uu1}
The set of trees\/ 
$\dU=\ens{t\app \uf ks}{k<\om\land s\in\bse\land t\in\bse}$
satisfies\/ \ref{ptf2} and\/ \ref{ptf3} while the union\/ 
$\dP\cup\dU$ is a\/ \sif.
\qed
\ele

\ble
\lam{uu2}
The set\/ $\dU$ is dense in\/ $\dU\cup\dP$. 
\ele
\bpf
Suppose that $T\in\dP$. 
The set $D(T)$ of all \mus s   
$\Phi=\sis{\vpi_k}{k\in\om}\in \spg\dP$,
such that $\vpi_k(\La)=T$ for some $k$, belongs to $\cM$ 
and obviously is dense in $\spg\dP$. 
It follows that $\Phi^j\in D(T)$ for some $j$, 
by the choice of $\dphi$. 
Then $\tx k\La=T$ for some $k$. 
However $\uf k\La\sq \tx k\La$. 
\epf

\ble
\lam{uu3}
If a set\/ $D\in\cM$, $D\sq\dP$ is pre-dense 
in\/ $\dP$, and\/ $U\in\dU$, then\/ $U\sqf\bigcup D$, 
that is, there is a finite\/ $D'\sq D$ with
$U\sq\bigcup D'$. 
Moreover\/ $D$ remains pre-dense in\/ $\dU\cup\dP$.
\ele
\bpf
Suppose that $U=\uf Ks\in\dU$, $K<\om$ and $s\in\bse.$ 
(The general case, when $U=t\app\uf Ks$ for some $t\in\bse,$ 
is easily redusible to the particular case $U=\uf Ks$ by 
substituting the set $t\app D$ for $D$.)
Consider the set $\Da\in\cM$ of all \mus s   
$\Phi=\sis{\vpi_k}{k\in\om}\in \spg\dP$ such that 
$K\in\abs\Phi$, $\lh s<h=\vys{\vpi_K}$, 
and for each $t\in 2^{h-1}$ there is a tree $S_t\in D$ with  
$\vpi_K(t)\sq S_t$.
The set $\Da$ 
is dense in $\spg\dP$ by Lemma~\ref{lsys} and 
the pre-density of $D$. 
Therefore there is an index $j$ such that $\Phi^j$ belongs 
to $\Da$.
Let this be witnessed by trees $S_t\in D\yt t\in 2^{h-1},$ 
where $\lh s<h=\vys{\vpj JK}$, so that $\vpj JK(t)\sq S_t$. 
Then 
$$
\textstyle
U=\uf Ks\sq\uf K\La\sq\bigcup_{t\in 2^{h-1}}\vpj JK(t)
\sq\bigcup_{t\in 2^{h-1}}S_t\sq\bigcup D'
$$ 
by construction, where $D'=\ens{S_t}{t\in 2^{h}}\sq D$ 
is finite.

To prove the pre-density, consider 
any string $t\in 2^{h-1}$ with $s\su t$.
Then $V=\uf Kt\in\dU$ and $V\sq U$. 
On the other hand, $V\sq S_t\in D$. 
Thus the tree $V$ witnesses that $U$ is compatible with 
$S_t\in D$ in $\dU\cup\dP$, as required.
\epf

\parf{Forcing a real away of a pre-dense set} 
\las{saway}

Let $\cM$ be still a countable transitive model of $\ZFC'$ 
and $\dP\in\cM\yd \dP\sq\pes$ be a \sif. 
The goal of the following Theorem~\ref{K} is to prove that, 
in the conditions of Definition~\ref{dPhi}, for any 
\dd\dP name $c$ of a real in $\dn,$ it is forced by the extended 
forcing $\dP\cup\dU$ that $c$ does not belong to sets $[U]$ where 
$u$ is a tree in $\dU$ --- unless $c$ is a name of one of 
reals in the \dd\Eo class of the generic real $x$ itself.
We begin with a suitable notation.

\bdf
\lam{rk}
A \rit{\dplo real name} is a system 
$\rc=\sis{\kc ni}{n<\om,\, i<2}$ of sets $\kc ni\sq\dP$ 
such that each set $C_n=\kc n0\cup \kc n1$ is 
dense or at least pre-dense in $\dP$ 
and if $S\in \kc n0$ and $T\in \kc n1$ then $S,T$ are 
incompatible in $\dP$.

If in addition $\sg\in\bse$ then define a \dplo real name  
$\sg\rc=\sis{\sg\app\kc ni}{n<\om,\, i<2}$, where 
$\sg\app\kc ni=\ens{\sg\app T}{T\in\kc ni}$. 

If a set $G\sq\dP$ is \dd\dP generic at least over 
the collection of all  sets $C_n$ then we define 
$\rc[G]\in\dn$ so that $\rc[G](n)=i$ iff $G\cap \kc ni\ne\pu$.
\edf

Thus any \dplo  real name $\rc=\sis{\kc ni}{}$ 
is a \dplo  name for a real in $\dn.$ 

Recall that $\dP$ adds a real $x\in\dn$.

\bpri
\lam{proj}
Let $k<\om$. 
Define a \dplo real name $\rpi=\sis{\kc ni}{n<\om\yi i<2}$
such that each set $\kc ni$ contains a single tree  
$\kR ni=\ens{s\in\bse}{\lh s>n\imp s(n)=i}\in\pes$.
Then $\rpi$ is a \dd\dP name of the \dd\dP generic real $x$, 
and accordingly each name $\sg\rpi$ 
($\sg\in\bse$) is a \dd\dP name of $\sg\app x$. 
\epri

Let $\rc=\sis{\kc ni}{}$ and $\rd=\sis{\kc ni}{}$ 
be \dplo real names. 
Say that  $T\in \pes$:
\bit
\item
\rit{directly forces\/ $\rc(n)=i$}, 
where $n<\om$ and $i=0,1$, iff $T\sq\kR ni$ 
(that is, the tree $T$ satisfies 
$x(n)=i$ for all $x\in[T]$); 

\item
\rit{directly forces\/ $s\su\rc$},  
where $s\in\bse,$ iff for all $n<\lh s$, $T$ 
directly forces $\rc(n)=i$, where $i=s(n)$; 

\item
\rit{directly forces\/ $\rd\ne\rc$}, iff there are strings 
$s,t\in\bse,$ incomparable in $\bse$ and such that  
$T$ directly forces $s\su\rc$ and $t\su\rd$; 

\item
\rit{directly forces\/ $\rc\nin[S]$},  
where $S\in\pet$, iff there is a string $s\in\bse\bez S$ 
such that $T$ directly forces $s\su \rc$; 
\vyk{
\item
\,{\ubf[applicable only for $\jta\in\plo$]}
\rit{weakly\/ \dd\plo forces\/ $\rc\nin[T]$},  
iff the set of all conditions 
$\jsg\in\plo$ that directly force\/ $\rc\nin[T]$ 
is dense in $\plo$ below $\jta$.
} 
\eit

\ble
\lam{K1}
If\/ $S_1,\dots,S_n,T\in \dP$ and\/ 
$\rc$ 
is a\/ \dplo real name then there exist trees\/ 
$S'_1,\dots,S'_n,T'\in\dP$ such that\/ $S'_i\sq S_i$ 
for all\/ $i=1,\dots,n$, $T'\sq T$, and\/ 
$T'$ directly forces\/ $\rc\nin[S']$, where\/ 
$S'=\bigcup_{1\le i\le n}S'_i$.
\ele
\bpf
Clearly there is a tree $T'\in\dP$, $T'\sq T$, which  
directly forces\/ $s\su\rc$ for some $s\in\bse$ 
satisfying $\lh s>\lh{(\roo {S_i})}$ for every $i$. 
Then there is a collection of strings $u_i\in S_i$ 
incomparable with $s$. 
Put $S'_i=S\ret {u_i}$; then obviously 
$s\nin S'=\bigcup_{1\le i\le n}S'_i$.
\epf

\ble
\lam{K2}
If\/ 
$\rc$ is a\/ \dplo real name, $\sg\in\bse,$
and\/ $T\in\dP$ directly forces\/ $\sg\rc\ne\rpi$, 
then there is a tree\/ 
$S\in\dP\yt S\sq T$, which directly forces\/ $\rc\nin[\sg\app S]$.
\ele
\bpf
Taking $T'=\sg\app T$ instead of $T$ and $\rc'=\sg\rc$ 
instead of $\rc$, we reduce the problem to the case 
$\sg=\La$, that is, $\sg\rc=\rc$ and $\sg\app S=S$.
Thus let's assume  that 
$T$ directly forces $\rc\ne\rpi$. 
There are incomparable strings $s,t\in\bse$ 
such that  
$T$ directly forces $s\su\rc$ and $t\su\rpi$.
Then by necessity $t\in T$, hence, $S=T\ret t\in \dP$ 
but $s\nin S$. 
By definition $S$ directly forces\/ $\rc\nin[S]$, 
as required.
\epf

\bte
\lam{K}
In the assumptions of Definition~\ref{dPhi}, suppose that\/ 
$\rc=\break
\sis{\kc mi}{m<\om\yi i<2}\in\cM$ is a\/ \dplo real name, 
and for every\/ $\sg\in\bse$ the set
$$
D(\sg)=\ens{T\in\dP}{T\text{ directly forces }\rc\ne{\sg\rpi}}
$$
is dense in\/ $\dP$. 
Let\/ $W\in\dP\cup\dU$ and\/ $U\in \dU$.
Then there is\/ 
a stronger condition\/ $V\in\dU\yd V\sq W$, which 
directly forces\/ $\rc\nin[U]$.
\ete
\bpf
By construction, $U=\sg\app\uf K{s_0}$, 
where $K<\om$ and $\sg,s_0\in\bse;$ 
we can assume that simply $s_0=\La$, so that $U=\sg\app\ufi K$. 
Further, by the same reasons as in the proof of Lemma~\ref{K2}, 
we can assume that  $\sg=\La$, so that $U=\ufi K$.
Finally, by Lemma~\ref{uu2}, we can assume that 
$W=\uf L{t_0}\in\dU$, where $L<\om$ and $t_0\in\bse.$   
The indices $K,L$ involved can be either equal or different.

There is an index $J$ such that the \mus\  
$\Phi^J=\sis{\vpj Jk}{k\in\om}$ satisfies 
$K,L\in\abs{\Phi^J}$ and 
$\vys{\vpj JL}>\lh{t_0}$, 
so that the trees 
$$
S_0=\vpj JK(\La)=\tx K{\La}\qand
T_0=\vpj JL(t_0)=\tx L{t_0} 
$$ 
in $\dP$ are defined. 
Note that $U\sq S_0$ and $W\sq T_0$.  

Consider the set $\cD$ of all \mus s 
$\Phi=\sis{\vpi_k}{k\in\om}\in \spg\dP$ such that 
$\Phi^J\cle\Phi$ and there is a tree 
$T\in\dP\yt T\sq T_0$
satisfying
\ben
\nenu
\atc\atc\atc\atc
\atc\atc\atc\atc
\itla{nen1}
$T$ directly forces $\rc\nin[S]$, where 
$S=\bigcup_{s\in 2^{h-1}}\vpi_K(s)$, 
$h=\vys{\vpi_K}$; \ and

\itla{nen2}
the tree $T$ occurs in $\Phi$ (see Section~\ref{spe}), and 
more specifically, 
$T=\vpi_L(t)$, 
where $t\in 2^{h'-1}$, $h'=\vys{\vpi_L}$,
and $t_0\su t$. 
\een

\ble
\lam{Ka}
$\cD$ is dense in\/ $\spg\dP$ above\/ $\Phi^J$. 
\ele
\bpf
Consider any \mus\ $\Phi^*=\sis{\vpi^*_k}{k\in\om}\in\spg\dP$ 
with $\Phi^J\cle\Phi^*$; the goal is to define a \mus\  
$\Phi'\in\cD$ such that $\Phi^*\cle\Phi'$.
By Corollary~\ref{csys} there is  an intermediate 
\mus\ $\Phi=\sis{\vpi_k}{k\in\om}\in\spg\dP$  
satisfying $\Phi^*\cll\Phi$; 
then any \mus\ $\Phi'\in\spg\dP$, which is a reduction of 
$\Phi$, still satisfies $\Phi^*\cll\Phi'$ and $\Phi^*\cle\Phi'$. 
Thus it suffices to find a \mus\ $\Phi'\in\cD$ 
which reduces $\Phi$.

Let $h=\vys{\vpi_K}$ and $h'=\vys{\vpi_L}$. 
Then $\vys{\vpj JK}<h$ and 
$\vys{\vpj JL}<h'$ strictly.
Pick a string $t\in 2^{h'-1}$ with $t_0\su t$;
let $R=\vpi_L(t)$; $R\sq T_0$ is a tree in $\dP$. 
Let $2^{h-1}=\ans{s_1,\dots,s_N}$, where $N$ is the 
integer $2^{h-1},$ 
and $S_i=\vpi_K(s_i)$. \vom

{\it Case 1\/}: $K\ne L$.
By Lemma~\ref{K1}, there exist trees  
$S'_1\sq S_1,\dots,S'_n\sq S_n$ and $T'\sq R$ in $\dP$ 
such that $T'$ directly forces\/ $\rc\nin[S']$, where 
$S'=\bigcup_{1\le i\le N}S'_i$.
Define a \mus\ $\Phi'=\sis{\vpi'_k}{k\in\om}\in\spg\dP$ so that 
$\vpi'_L(t) = T'$, $\vpi'_K(s_i)=S'_i$ for all $i=1,...,N$, 
and $\vpi'_k(s)=\vpi_k(s)$ for all other applicable values 
of $k$ and $s$. 
Then $\Phi'$ belongs to $\cD$ and is a reduction 
of $\Phi$, as required. 
\vom

{\it Case 2\/}: $L=K$, and hence $h'=h$.
Now $t$ is one of $s_i$, say $t=s_{i(t)}$, and the construction 
as in Case 1 does not work. 
Nevertheless, following the same arguments, we find trees  
$S'_i\sq S_i$, $i=1,\dots,N$, $i\ne i(t)$, 
and $T\sq R=\vpi_K(t)$ in $\dP$ 
such that $T$ directly forces\/ $\rc\nin[S']$, where 
$S'=\bigcup_{1\le i\le N,\,i\ne i(t)}S'_i$.

Further, as the set $D(\La)$ is dense, there is a tree $T'\in\dP$, 
$T'\sq T$, which directly forces $\rc\ne{\rpi}$. 
By Lemma~\ref{K2}, there is an even smaller tree $T''\in\dP$, 
$T''\sq T'$, which directly forces $\rc\nin[T'']$, that is, 
$T''$ directly forces $\rc\nin[S'\cup T'']$. 
Define a \mus\ $\Phi'=\sis{\vpi'_k}{k\in\om}\in\spg\dP$ so that 
$\vpi'_K(s_i)=S'_i$ for all $i=1,...,N$, $i\ne i(t)$, 
$\vpi'_K(t)=T''$,  
and $\vpi'_k(s)=\vpi_k(s)$ for all other applicable values 
of $k$ and $s$. 
Then $\Phi'\in\cD$ and $\Phi'$ is a reduction 
of $\Phi$. 
\epF{Lemma}

Come back to the proof of the theorem.
It follows from the lemma that there is an index $j\ge J$ 
such that the \mus\ $\Phi^j=\sis{\vpj jk}{k\in\om}$ 
belongs to $\cD$, and let this be witnessed by a tree 
$T=\vpj jL(t)\sq T_0=\vpj JL(t_0)=\tx L{t_0}$, 
where $t\in 2^{h'-1}$, $h'=\vys{\vpj jL}$,
and $t_0\su t$, satisfying \ref{nen1}. 

Consider the tree $V=\uf L{t}\in\dU$. 
By construction we have both $V\sq W$ and $V\sq T\sq T_0$.
Therefore $V$ directly forces $\rc\nin[S]$ by the choice 
of $T$ (which satisfies \ref{nen1}), where 
$S=\bigcup_{s\in 2^{h-1}}\vpj JK(s)$, $h=\vys{\vpi_K}$.
And finally, we have $U\sq S$, so that 
$V$ directly forces $\rc\nin[S]$, as required. 
\epf

\parf{Jensen's forcing} 
\las{jfor}

In this section, 
{\ubf we argue in $\rL$, the constructible universe.}
Let $\lel$ be the canonical wellordering of $\rL$.

\bdf
[in $\rL$]
\lam{uxi}
Following the construction in
\cite[Section 3]{jenmin} \rit{mutatis mutandis}, 
we define, by induction on $\xi<\omi$, a countable set  of trees 
$\dU_\xi\sq\pes$ satisfying requirements \ref{ptf2} and 
\ref{ptf3} of Section \ref{tre}, 
as follows.

Let $\dU_0$ consist of all clopen trees $\pu\ne S\sq\bse$, 
including $\bse$ itself.

Suppose that $0<\la<\omi$, and countable sets 
$\dU_\xi\sq\pes$ are already defined. 
Let $\cM_\xi$ be the least model $\cM$ of $\zfc'$ of the form 
$\rL_\ka\yt\ka<\omi$, 
containing $\sis{\dU_\xi}{\xi<\la}$ and such that 
$\al<\omi^\cM$ and all sets $\dU_\xi$, $\xi<\la$, 
are countable in $\cM$.

Then $\dP_\la=\bigcup_{\xi<\la}\dU_\xi$ is countable in $\cM$, too. 
Let $\sis{\Phi^j}{j<\om}$ be the $\le_\rL$-least sequence of 
\mus s $\Phi^j\in\spg{\dP_\la}$, \dd\cle increasing and generic 
over $\cM_\la$, and let $\dU_\la=\dU$ be defined, 
as in Definition~\ref{dPhi} and Lemma~\ref{uu1}. 

Let  $\dP=\bigcup_{\xi<\omi}\dU_\xi$. 
\edf

\bpro
[in $\rL$]
\lam{uxip}
The sequence\/ $\sis{\dU_\xi}{\xi<\omi}$ belongs to\/ 
$\id{\hc}1$.\qed
\epro

\ble
[in $\rL$]
\lam{jden}
If a set\/ $D\in\cM_\xi\yt D\sq {\dP_\xi}$ is 
pre-dense in\/ ${\dP_\xi}$ then it remains pre-dense in\/ 
$\dP$. 
Hence if\/ $\xi<\omi$ then\/ 
${\dU_\xi}$ is pre-dense in\/ $\dP$.
\ele
\bpf
By induction on $\la\ge \xi$, 
if $D$ is pre-dense in ${\dP_\la}$ then it 
remains pre-dense in 
${\dP_{\la+1}}=\dP_\la\cup\dU_\la$ 
by Lemma~\ref{uu3}. 
Limit steps are obvious. 
To prove the second part, note that 
${\dU_\xi}$ is dense in ${\dP_{\xi+1}}$ by Lemma~\ref{uu2},
and $\dU_\xi\in\cM_{\xi+1}$.
\epf

\ble
[in $\rL$]
\lam{club}
If\/ $X\sq\HC=\rL_{\omi}$ then the set\/ $W_X$ of all 
ordinals\/ $\xi<\omi$ such that\/ 
$\stk{\rL_\xi}{X\cap\rL_\xi}$ is an elementary submodel of\/  
$\stk{\rL_{\omi}}{X}$ and\/ $X\cap\rL_\xi\in\cM_\xi$ 
is unbounded in\/ $\omi$.
More generally, if\/ $X_n\sq\HC$ for all\/ $n$ 
then the set\/ $W$ of all 
ordinals\/ $\xi<\omi$, such that\/ 
$\stk{\rL_\xi}{\sis{X_n\cap\rL_\xi}{n<\om}}$ 
is an elementary submodel of\/  
$\stk{\rL_{\omi}}{\sis{X_n}{n<\om}}$ 
and\/ $\sis{X_n\cap\rL_\xi}{n<\om}\in\cM_\xi$, 
is unbounded in\/ $\omi$.
\ele
\bpf
Let $\xi_0<\omi$. 
By standard arguments, there are ordinals $\xi<\la<\omi$, 
$\xi>\xi_0$, such 
that $\stk{\rL_\la}{\rL_\xi,X\cap\rL_\xi}$ is an elementary 
submodel of $\stk{\rL_{\om_2}}{\rL_{\omi},X}$.
Then $\stk{\rL_\xi}{X\cap\rL_\xi}$ is an elementary submodel of   
$\stk{\rL_{\omi}}{X}$, of course. 
Moreover, $\xi$ is uncountable in $\rL_\la$, hence 
$\rL_\la\sq\cM_\xi$. 
It follows that $X\cap\rL_\xi\in\cM_\xi$ since 
$X\cap\rL_\xi\in\rL_\la$ by construction.
The second claim does not differ much.
\epf

\bcor
[compare to \cite{jenmin}, Lemma 6]
\lam{ccc}
The forcing\/ $\dP$ satisfies CCC in $\rL$.
\ecor
\bpf
Suppose that $A\sq\dP$ is a maximal antichain. 
By Lemma~\ref{club}, there is an ordinal $\xi$ such that 
$A'=A\cap{\dP_\xi}$ is a maximal antichain in ${\dP_\xi}$ 
and $A'\in\cM_\xi$. 
But then $A'$ remains pre-dense, therefore, 
still a maximal antichain, in the 
whole set $\dP$ by Lemma~\ref{jden}. 
It follows that $A=A'$ is countable.
\epf

\parf{The model} 
\las{mod}

We consider the set $\dP\in\rL$ (Definition~\ref{uxi})   
as a forcing notion over $\rL$. 
 
\ble
[compare to Lemma 7 in \cite{jenmin}]
\lam{mod1}
A real\/ $x\in\dn$ is\/ $\dP$-generic over\/ $\rL$ iff\/ 
$x\in Z=\bigcap_{\xi<\omi^\rL}\bigcup_{U\in\dU_\xi}[U]$. 
\ele
\bpf
If $\xi<\omi^\rL$ then $\dU_\xi$ is pre-dense in $\dP$ 
by Lemma~\ref{jden}, therefore any real $x\in\dn$ 
$\dP$-generic over $\rL$ belongs to $\bigcup_{U\in\dU_\xi}[U]$. 

To prove the converse, suppose that $x\in Z$ and prove that 
$x$ is $\dP$-generic over $\rL$. 
Consider a maximal antichain $A\sq\dP$ in $\rL$; we have to 
prove that $x\in \bigcup_{T\in A}[T]$. 
Note that $A\sq\dP_\xi$ for some 
$\xi<\omi^\rL$ by Corollary~\ref{ccc}.
But then every tree $U\in\dU_\xi$ satisfies 
$U\sqf\bigcup A$ by Lemma~\ref{uu3}, so that 
$\bigcup_{U\in \dU_\xi}[U]\sq \bigcup_{T\in A}[T]$, and hence 
$x\in \bigcup_{T\in A}[T]$, as required.
\epf       

\bcor
[compare to Corollary 9 in \cite{jenmin}]
\lam{mod2}
In any generic extension of\/ $\rL$, the set of all reals 
in\/ $\dn$ $\dP$-generic over\/ $\rL$ is\/ $\ip\HC1$ and\/ 
$\ip12$. 
\ecor
\bpf
Use Lemma~\ref{mod1} and Proposition~\ref{uxip}.
\epf  

\bdf
\lam{gg}
From now on, let $G\sq\dP$ be a set \dd\dP generic over $\rL$,
so that $X=\bigcap_{T\in G}[T]$ is 
a singleton $X_G=\ans{x_G}$.
\edf

Compare the next lemma to Lemma 10 in \cite{jenmin}. 
While Jensen's forcing notion in \cite{jenmin} guarantees 
that there is a single generic real in the extension, 
the forcing notion $\dP$ we use adds a whole \dd\Eo class 
(a countable set) of generic reals!

\ble
[in the assumptions of Definition~\ref{gg}] 
\lam{only}
If\/ $y\in\rL[G]\cap\dn$ then\/ $y$ is a\/ $\dP$-generic real 
over\/ $\rL$ iff\/ $y\in\eko{x_G}=\ens{\sg\app x_G}{\sg\in\bse}$.
\ele
\bpf
The real $x_G$ itself is \dd\dP generic, of course. 
It follows that any real $y=\sg\app x_G\in\eko{x_G}$ 
is \dd\dP generic as well 
since the forcing $\dP$ is 
by definition invariant under the action of any $\sg\in\bse.$

To prove the converse, suppose towards the contrary that 
there is a tree $T\in\dP$ and 
a \dd\dP real name $\rc=\sis{\kc ni}{n<\om,\,i=0,1}\in\rL$ 
such that $T$ \dd\dP forces that $\rc$ is \dd\dP generic 
while $\dP$ forces that $\rc\ne\sg\app\rpi$ for all 
$\sg\in\bse.$  

Let $C_n=\kc n0\cup\kc n1$; this is a pre-dense set in $\dP$. 
It follows from Lemma~\ref{club} that there is an ordinal 
$\la<\omi$ such that 
each set 
$C'_n=C_n\cap {\dP_\la}$ is pre-dense in ${\dP_\la}$, 
and the sequence $\sis{\xc ni}{n<\om,\,i=0,1}$ belongs to 
$\cM_\la$, where $\xc ni=C'_n\cap \kc ni$ --- 
then $C'_n$ is pre-dense in $\dP$ too, by Lemma~\ref{jden}. 
Thus we can assume that in fact $C_n=C'_n$, that is, 
$\rc\in\cM_\la$ and $\rc$ is a \dd{\dP}real name.

Further, as $\plo$ forces that $\rc\ne\sg\app\rpi$, 
the set $D(\sg)$ 
of all conditions $S\in\dP$ which directly force 
$\rc\ne\sg\app\rpi$, is dense in $\dP$ --- 
for every $\sg\in\bse.$  
Therefore, still by Lemma \ref{club}, 
we may 
assume that the same ordinal $\la$ as above satisfies the 
following: each set $D'(\sg)=D(\sg)\cap{\dP_\la}$ is dense in 
${\dP_\la}$. 

Applying Theorem~\ref{K} with $\dP=\dP_\la$, $\dU=\dU_\la$, 
and $\dP\cup\dU=\dP_{\la+1}$, we conclude that for each 
$U\in\dU_\la$ the set $Q_U$ of all conditions 
$V\in \dP_{\la+1}$ which directly force $\rc\nin[U]$, 
is dense in $\dP_{\la+1}$.
As obviously $Q_U\in\cM_{\la+1}$, we further conclude that 
$Q_U$ is pre-dense in the whole forcing $\dP$ 
by Lemma~\ref{jden}.
This implies that $\dP$ forces 
$\rc\nin\bigcup_{U\in \dU_\la}[U]$, 
hence, forces that $\rc$ is not \dd\dP generic, 
by Lemma~\ref{mod1}.
But this contradicts to the choice of $T$.
\epf


\ble
[in the assumptions of Definition~\ref{gg}] 
\lam{sym}
$x_G$ is not\/ $\od$ in\/ $\rL[G]$. 
\ele
\bpf
Suppose towards the contrary that 
there is a tree $T\in G$ and a formula $\vt(x)$ with 
ordinal parameters such that $T$ \dd\dP forces that $x_G$ 
is the only $x\in\dn$ satisfying $\vt(x)$. 
Let $s=\roo T$, so that both $s\we0$ and $s\we1$ belong to $T$.
Then either $s\we0\su x_G$ or $s\we1\su x_G$; let, say, 
$s\we0\su x_G$. 

Let $n=\lh s$ and $\sg=0^n\we 1$, so that all three strings 
$s\we0,s\we1,\sg$ belong to $2^{n+1},$ and $s\we1=\sg\app s\we0$.
As the forcing $\dP$ is invariant under the action of $\sg$, 
the set $G'=\sg\app G$ is \dd\dP generic over $\rL$, and 
$T=\sg\app T\in G'$. 
It follows that it is true in $\rL[G']=\rL[G]$ that the real 
$x'=x_{G'}=\sg\app x_G$ is still the only $x$ satisfying $\vt(x)$.
However obviously $x'\ne x$!
\epf

Now, arguing in the \dd\dP generic model 
$\rL[G]=\rL[x_G]$, 
we observe that the countable set $X=\eko{x_G}$ 
is exactly the set of 
all \dd\dP generic reals by Lemma~\ref{only}, hence it 
belongs to $\ip12$ by Corollary~\ref {mod2}, and finally 
it contains no $\od$ elements by Lemma~\ref{sym}, as 
required.\vom

\qeD{Theorem~\ref{mt}} 

\back
The authors thank Jindra Zapletal and Ali Enayat for 
fruitful discussions.
\eack

\bibliographystyle{plain}

\end{document}